\documentclass[11pt]{article}\hfuzz=10pt \topmargin=-0.5cm\hfuzz=10pt  \oddsidemargin=-0.2cm \textheight 229mm \textwidth=16.5cm
\newcommand{\comm}[1]{}

\usepackage{amssymb,epsfig}
\usepackage{amsthm}
\usepackage[numbers]{natbib}
\def\citet{\cite}

    \csname @addtoreset\endcsname{figure}{section}
    \csname @addtoreset\endcsname{table}{section}

\newtheorem{theorem}{Theorem}[section]
\newtheorem{lemma}{Lemma}[section]

\newtheorem{corollary}{Corollary}[section]

\newtheorem{definition}{Definition}[section]
\newtheorem{remark}{Remark}[section]

\def\e{\varepsilon}

\def\defi{\stackrel{{\scriptscriptstyle \Delta}}{=}}

\def\o{\omega}
\def\O{\Omega}

\def\Y{{\cal Y}}

\def\w{\widehat}
\def\Ind{{\,\rm Ind\,}}
\def\Ind{{\mathbb{I}}}

\def\R{{\bf R}}
\def\E{{\bf E}}

\def\Z{{\cal Z}}

\def\H{{\cal H}}

\def\b{\beta}
\def\s{\delta}

\def\C{{\bf C}}

\def\ww{\widetilde}
\def\X{{\cal X}}
\def\t{\theta}
\def\oo{\bar}
\def\s{\sigma}

\newcommand{\be}{\begin{equation}}
\newcommand{\ee}{\end{equation}}
\newcommand{\bd}{\begin{displaymath}}
\newcommand{\ed}{\end{displaymath}}
\newcommand{\ba}{\begin{array}{ll}}
\newcommand{\ea}{\end{array}}
\newcommand{\baa}{\begin{eqnarray}}
\newcommand{\eaa}{\end{eqnarray}}
\newcommand{\baaa}{\begin{eqnarray*}}
\newcommand{\eaaa}{\end{eqnarray*}}
\font\sm=cmr10

\def\oo{\bar}

\def\Q{{\cal Q}}
\def\Re{{\rm Re\,}}

\def\sinc{{\rm sinc\,}}
\def\ew{\left(e^{i\o}\right)}
\def\T{{\mathbb{T}}}
\def\ZZ{{\mathbb{Z}}}

\def\TT{{\cal T}}
\def\BL{{{\scriptscriptstyle BL}}}
\def\lBL{{\scriptscriptstyle LBL}}

\def\BS{{\scriptscriptstyle BS}}

\def\spl{{\scriptscriptstyle spline}}
\def\MA{{\scriptscriptstyle MA}}
\def\XN{B_N}
\def\BN{{\cal B}_N}
\def\BNN{{\cal B}_\infty}
\def\BN{{\mathbb{B}}_N}
\def\BNN{{\mathbb{B}}_\infty}
\def\LBL{\ell_2^\lBL(-\infty,s)}
\def\LBLN{\ell_{2,N}^\lBL(-\infty,s)}
\title{
On  causal extrapolation of sequences with  applications to forecasting\footnote{Accepted for publication in 
Applied Mathematics and Computation}}\author{
Nikolai Dokuchaev\footnote{
 Department of Mathematics \& Statistics, Curtin
University, GPO Box U1987, Perth, 6845 Western Australia\comm{,  and National Research University ITMO, 197101 Russia.  Email: N.Dokuchaev@curtin.edu.au}}}

\begin{document}
\def\break{}%
\def\brea{}
\def\breakk{}
\maketitle
\begin{abstract}
  The paper suggests a method of extrapolation of notion of one-sided semi-infinite
sequences representing traces of two-sided band-limited sequences; this features ensure uniqueness
of this extrapolation and possibility to use this for forecasting.
This lead to a forecasting method
for more general  sequences without this feature based on minimization of the mean square error
between the observed path and a predicable sequence.
These procedure involves calculation of this predictable path; the procedure can be interpreted as causal smoothing. The corresponding smoothed
sequences allow  unique
extrapolations to future times that can be interpreted as optimal forecasts.
\par
 {\bf Key words}:  smooth sequences,
extrapolation, frequency analysis, forecasting
\par
AMS 2010 classification:
42A38, 
93E10, 
	42A99 

\end{abstract}
\section{Introduction}  We study causal dynamic approximation and extrapolation
of real sequences in deterministic setting, i.e.  without probabilistic assumptions.
 Extrapolation of sequences can be used for forecasting and was
studied intensively, for example,  in the framework of system identification methods; see e.g. \cite{SFS}.
In signal processing, there is an approach oriented on the
frequency analysis and exploring special features of the band-limited processes such as predictability. For stochastic
 stationary discrete time processes, the connection between  predicability and degeneracy of the spectrum  was established by
the classical Szeg\"o-Kolmogorov Theorem;  see a  recent reviews in \citet{Bin}. This theorem says that
the optimal prediction error is zero  if its spectral density is vanishing with a certain rate at a point of the
unit circle $\T=\{z\in\C:\ |z|=1\}$, in particular, if it is vanishing on an arc on $\T$.
In this case,
the process is called
"band-limited". This result was expanded on more
general stochastic processes featuring
spectral densities; see, e.g., \citet{CS}.
In deterministic pathwise setting without probability assumptions, this result was expanded  \citet{D10,D12a} on sequences with
 Z-transform vanishing at a point of $\T$.

The present suggests to use the predicability featuring by band-limited processes
for forecasting of more general processes that are not necessarily band-limited. This requires  calculation of a
trace of a band-limited process representing optimal approximation of the
available observations.  The extrapolation of this  trace of a band-limited process can be used as a forecast.
The motivation for that approach is based on the assumption that a band-limited part of a process can be
interpreted as its regular part purified from a noise represented by high-frequency component.
This leads to a problem of causal band-limited approximations for non-bandlimited
processes which can be interpreted as a  causal band-limited smoothing.

 A known two-sided  sequence can be converted into a band-limited process with a low-pass filter, and the resulting process will
be an optimal band-limited approximation.     However,
 a ideal low-pass filter is non-causal; therefore, it cannot be applied for  dynamically  observable processes with
  unavailable  future values which excludes  predicting and extrapolation problems.
Respectively,  causal smoothing cannot convert a process into a band-limited one; it is known that the distance of an ideal low-pass
filter from the set of all causal filters is positive \cite{rema}.
 There are many works devoted to causal smoothing and sampling, oriented on estimation and minimization of errors in $L_2$-norms or similar norms, especially in stochastic setting; see e.g.\citet{AU,Alem,rema,CTao,CJR,D10,D12a,D16,PFG,jerri,K,W,Zhao, Zhao2}.

 The present paper readdresses  the problem of causal  band-limited smoothing
   approximation and considers the problem of causal  band-limited extrapolation for one-sided real sequences that are not necessary
paths of band-limited processes.

We  consider purely discrete time processes rather than samples of continuous time processes;
one may say that the values between fixed discrete times are  are not included into consideration.
This setting imposes certain restrictions. In particular, it does not allow to consider continuously
variable locations of the sampling points, as is common in sampling analysis of continuous time processes;
see e.g. \citet{BS,PFG,FKR,LeeF}.  For continuous
time processes, the predicting horizon can be selected to be arbitrarily  small, such as in the model considered  in \citet{BS};
this possibility is absent for discrete time processes considered below. In addition,
 it is not obvious how to
define for discrete time processes or sequences an analog of the continuous time analyticity that is often associated with
predicability.

Further, we consider the problem  in the deterministic setting, i.e. pathwise.
This means that the method has to rely on the intrinsic properties of a sole underlying sequence without appealing
to statistical properties of an ensemble  of sequences. In particular, we use a
pathwise optimality criterion rather than criterions calculated via the expectation on a probability space such as
mean variance criterions.

In addition, we consider an approximation that does not target the match of the values at any set of selected points; the error is not expected to be small.
This is different from  a more common setting where the goal is to match
an approximating curve with the underlying process
at certain sampling points; see e.g. \citet{CJR,FKR,jerri,LeeF}.
Our setting is closer to the setting from \citet{PFG,PFG1,TH,Zhao,Zhao2}.
In \citet{PFG,PFG1}, the point-wise matching error
was estimated for a sampling series and for a band-limited process representing
smoothed
underlying  continuous time process; the  estimate featured a
given vanishing error. In \citet{TH}, the problem  of minimization of the total energy  of the approximating bandlimited process was considered; this causal approximation was constructed
 within  a given distance from  the original process smoothed by an ideal low-pass filter.
 Another related result was obtained in \citet{F94}, where an interpolation problem for absent sampling points  was considered in a setting with vanishing error, for a finite number of sampling points. In \citet{Zhao},  extrapolation of a trace of a band-limited process from a finite number of points was considered in a frequency
 setting for a general linear transform and some special Slepian's type basis in the frequency domain. In \citet{Zhao2}, a setting  similar to \citet{Zhao}
 was considered for extrapolation of a trace of continuous time process from a finite interval using a special basis from eigenfunctions  in the frequency domain.
 Our setting is different:
we consider extrapolation without exact match of values for the underlying process. Therefore, we suggest
to calculate extrapolations that can be used for forecasting that
 are not necessarily
matching the values of the underlying process. This allows to consider
semi-infinite underlying processes that are not paths of band-limited processes.

\par

It can be noted that the framework of two-sided sequences
 required for
detecting of the bandlimitness via Z-transforms are  not always convenient to use. For example, consider a situation where the data is collected
dynamically during a prolonged time interval. For many models, it is more convenient
 to represent this data  flow
as one-sided sequences such that $x(t)$  represents outdated observations with diminishing significance as
$t\to -\infty$. However, application of the two-sided Z-transform requires to select some
past time at the middle of the time interval of the observations as the zero point for a model of  the two-sided sequence; this could be inconvenient.
On the other hand, a straightforward application of the one-sided
Z-transform to the historical data represented as
 one-sided sequences  generates Z-transforms that cannot vanish on a part of the unit circle, even for
 traces of band-limited two-sided sequences.
So far,  the notion of bandlimitness was not expanded  on the one-sided sequences  $\{x(t)\}_{t=0,-1,-2,...,-\infty}$.
  The paper addresses this  problem, as well as the problem of detecting
one-sided semi-infinite sequences that can be extended into two-sided band-limited processes (i.e. representing traces of band-limited processes).

The paper suggests a method of extrapolation of notion of one-sided semi-infinite
sequences representing traces of two-sided band-limited sequences; this features ensure uniqueness
of this extrapolation and possibility to use this for forecasting.
This lead to a forecasting method
for more general  sequences without this feature based on minimization of the mean square error
between the observed path and a predicable sequence.
These procedure involves calculation of this predictable path; the procedure can be interpreted as causal smoothing. The corresponding smoothed
sequences allow  unique
extrapolations to future times that can be interpreted as optimal forecasts.

For the solution, we use  non-singularity  of special sinc matrices  obtained in \citet{LeeF}
for the solution of the so-called superoscillations problem for continuous time processes; see the references in  \citet{FKR,LeeF}. It can be noted that the setting in \citet{FKR,LeeF}
considers exact matching of the band-limited process and the underlying process in certain points, which is  different from our setting.

 The sustainability of the  method is demonstrated with  some  numerical experiments where we compare the band-limited extrapolation with some classical spline based interpolations.
\section{Definitions}
We use notation $\sinc(x)=\sin(x)/x$, and we denote by $\ZZ$  the set of all integers. \par
We assume that we are given $\O\in(0,\pi)$ an a positive integer $N$. In addition, we are given $s\in\ZZ$ and $q\in\{k\in\ZZ:\  k<s\}\cup\{-\infty\}$.
\par
Let $\TT=\{t\in\ZZ: \ q\le t\le s\}$ if $q>-\infty$ and
$\TT=\{t\in\ZZ: \ t\le s\}$ if $q=-\infty$.
\par
Let $\ZZ_N$ be the set of all integers $k$ such that $|k|\le N$.
\par
 For a Hilbert space $H$, we denote by $(\cdot,\cdot)_{H}$ the
corresponding inner product. We denote by $L_2(D)$ the usual Hilbert space of complex valued
square integrable functions $x:D\to\C$, where $D$ is an interval  in $\R$.   We denote by $\ell_r$ the set of all
sequences $x=\{x(t)\}_{t\in\ZZ}\subset\C$, such that
$\|x\|_{\ell_r}=\left(\sum_{t=-\infty}^{\infty}|x(t)|^r\right)^{1/r}<+\infty$
for $r\in[1,\infty)$ or  $\|x\|_{\ell_\infty}=\sup_t|x(t)|<+\infty$
for $r=+\infty$. \par Let $\ell_r^+$ be the set of all sequences
$x\in\ell_r$ such that $x(t)=0$ for $t=-1,-2,-3,...$; see, e.g. \citet{Yosida}.

\par Let $\T=\{z\in\C:\  |z|=1\}$.
\par
For  $x\in \ell_1$ or $x\in \ell_2$, we denote by $X=\Z x$ the
Z-transform  \baaa X(z)=\sum_{t=-\infty}^{\infty}x(t)z^{-t},\quad
z\in\T. \eaaa Respectively, the inverse Z-transform  $x=\Z^{-1}X$ is
defined as \baaa x(t)=\frac{1}{2\pi}\int_{-\pi}^\pi
X\left(e^{i\o}\right) e^{i\o t}d\o, \quad t=0,\pm 1,\pm 2,....\eaaa
If $x\in \ell_2$, then $X|_\T$ is defined as an element of
$L_2(\T)$.

Let $\tau\in\ZZ\cup\{+\infty\}$ and $\t<\tau$; the case where
$\t=-\infty$  is not excluded. We denote by $\ell_2(\t,\tau)$ the
Hilbert space of complex valued sequences $\{x(t)\}_{t=\t}^\tau$
such that
$\|x\|_{\ell_2(\t,\tau)}=\left(\sum_{t=\t}^\tau|x(t)|^2\right)^{1/2}<+\infty$.
\par
Let $\Y_N$ be the Hilbert space of sequences $\{y_k\}_{k=
-N}^{N}\subset\C$ provided with the $\ell_2$-norm, i.e.,
$\|y\|_{\Y_N}=\left(\sum_{k\in \ZZ_N }|y_k|^2\right)^{1/2}<+\infty$.

\par
Let $\ell_2^\BL$ be the set of all mappings $x\in\ell_2$  such
that $X\ew \in L_2(-\pi,\pi)$ and $X\ew =0$ for $|\o|>\O$, where $X=\Z x$. We will call the the corresponding processes $x=\Z^{-1}X$
{\em band-limited}.

\par Let $\BN$ be  the set of all $X\in
L_2(\TT)$ such that there exists a sequence
$\{y_k\}_{k=-N}^N\in\Y_N$ such that $X\ew
=\sum_{k=-N}^{N}y_ke^{ik\o\pi/\O}\Ind_{\{|\o|\le\O\}}$, where $\Ind$ is
the indicator function.

Consider the Hilbert spaces of sequences $\X=\ell_2$ and
$\X_-=\ell_2(q,s)$.

 Let $\XN$ be the subset of $\X_-$ consisting of sequences
$\{x(t)\}_{t\in\TT}$, where $x\in\X$ are such that $x(t)=(\Z^{-1}
X)(t)$ for $t\in\TT$ for some $X\ew \in \BN $.

\begin{definition}\label{defLBL}
We call an one-sided sequence $x'\in\ell_2(-\infty,s)$  left band-limited if
there exists  $x\in\ell_2^\BL$  such that $x(t)=x'(t)$ for $t\le s$.
We denote by $\LBL$ the set of all these sequences $x'$, and we denote
by $\LBLN$  the set of all these sequences $x'$ such that the corresponding extrapolation  $x$ belongs to $\XN$.
\end{definition}

\section{Main results}
\subsection{Uniqueness of the extrapolation for left band-limited processes}
\begin{lemma}\label{propU} \begin{itemize}\item[(i)] For $x_L\in\LBL$,  the extrapolation $\oo x\in\ell_2$ described in Definition
    \ref{defLBL} is uniquely defined.
  \item[(ii)]
  If $s-q\ge 2N+1$, then,  for
  any $x\in\XN $, there exists an
unique  $X\in\BN$ such that $x(t)=(\Z^{-1} X)(t)$ for $t\in\TT$. \end{itemize}
\end{lemma}
\par
By Lemma \ref{propU} (i), the future  of a
band-limited process $x$
 is uniquely defined by its  history
$\{x(t),\ t\le s\}$.  This statement represent a reformulation in the deterministic setting
of  the classical Szeg\"o-Kolmogorov Theorem for stationary Gaussian processes
\citet{Bin}. In addition, Lemma \ref{propU}(ii) states that the future of  processes
from $\XN $ is uniquely defined by a finite set of historical
values that has at least $2N+1$ elements.
\begin{corollary}\label{rem1}
\begin{itemize}
\item[(i)]  If $q>-\infty$ and  $s-q\le  2N+1$, then
 $\{x(t)\}_{t\in\TT}\in \XN $ for any $x\in\ell_2$.
 \item[(ii)] If   $s-q< 2N+1$, then there are many
$X\in\BN$ such that $x(t)=(\Z^{-1} X)(t)$  for $t\in\TT$; they form a linear
manifold in $\XN $.
\end{itemize}
\end{corollary}
\subsection{Optimal
left band-limited approximation} Let $x\in\X$ be a process. We assume
that the sequence $\{x(t)\}_{t\in\TT}$ represents available
historical data.
\begin{theorem}\label{Th1}\begin{itemize}
\item[(i)] There exists an optimal solution  $\w x$
of the minimization problem \baa &&\hbox{Minimize}\quad  \sum_{t=q}^s|\w
x(t)-x(t)|^2  \quad\breakk\hbox{over}\quad \w x\in \XN .\label{min} \eaa
\item[(ii)]
If  $s-q\ge
2N+1$, then  the corresponding optimal process $\w x$ is uniquely
defined.
\end{itemize}
\end{theorem}
\par
Under the assumptions of Theorem \ref{Th1}(ii), by  Lemma \ref{propU},  there exists a unique extrapolation of the
band-limited solution $\w x$ of problem (\ref{min}) on the future
times  $t>s$. It can
be interpreted  as the optimal forecast (optimal given $\O$ and
$N$).
\begin{corollary}\label{rem2}
If $s-q<2N+1$ then there are
many optimal processes $\w x\in \XN$ such that $ \sum_{t=q}^s|\w
x(t)-x(t)|^2=0$; they  form a linear manifold in
$\XN $.
\end{corollary}
Up to the end of  this section, we assume that the assumptions of Lemma \ref{propU} and Theorem \ref{Th1}(ii) are satisfied.
\subsubsection*{The optimal solution}
Let  the operator $\Q: \Y_N\to \XN $ be defined as $\w x=\Q y=\Z^{-1}\w X$, where
\baaa
\w X\ew =\sum_{k\in
\ZZ_N }y_ke^{ik\o\pi/\O}\Ind_{\{|\o|\le\O\}},
\label{wX}\eaaa  for the corresponding  $y=\{y_k\}\in \Y_N$.
Similarly to the classical sinc representation, we obtain that \baa \w x(t)=\frac{1}{2\pi}
\int_{-\O}^{\O}\left(\sum_{k\in \ZZ_N }y_k e^{ik\o\pi/\O}\right)e^{i\o
t}d\o\brea=\frac{1}{2\pi}
\sum_{k\in \ZZ_N }y_k\int_{-\O}^{\O}e^{ik\o\pi/\O+i\o t}d\o\nonumber\\
=\frac{1}{2\pi}\sum_{k\in \ZZ_N }y_k \frac{e^{ik\pi+i\O t}-
e^{-ik\pi-i\O t}}{ik\pi/\O+it}\brea=\frac{\O}{\pi}\sum_{k\in
\ZZ_N }y_k \sinc(k\pi+\O t)=(\Q y)(t).\label{sinc}\eaa
It follows that the $\Q: \Y_N\to \XN $ is actually defined as
\baaa \w x(t)=(\Q y)(t)=\frac{\O}{\pi}\sum_{k\in
\ZZ_N }y_k \sinc(k\pi+\O t).\label{Qs}\eaaa  Consider the operator  $\Q^*:\XN \to \Y_N$  being adjoint to the operator
$\Q:\Y_N\to\XN $, i.e., such that
\baa
(\Q^*x)_k=\frac{\O}{\pi}\sum_{t\in\TT}\sinc(k\pi+\O t)x(t).
\label{Q*}\eaa
By the property of the sinc function, it follows that this convolution maps continuously $\ell_2$ into $\ell_2$. Hence
the operator $\Q^*$ can be extended as a continuous linear operator $\Q^*:\X_-\to \Y_N$.

Consider the
linear bounded non-negative definite  Hermitian operator $R:\Y_N\to \Y_N$  defined as
\baaa
R=\Q^*\Q.
\eaaa
\begin{theorem}
\label{ThP}
\begin{itemize}\item[(i)] The operator $R:\Y_N\to\Y_N$ has a bounded inverse
 operator $R^{-1}:\Y_N\to\Y_N$.
\item[(ii)] Problem (\ref{min}) has a unique solution \baa
\w x=\Q R^{-1} \Q^*x.\label{wx}
\eaa
\end{itemize}
\end{theorem}
\begin{remark}
It can be noted that $\w x=\Q \Q^+ x$, where $\Q^+=R^{-1} \Q^*:\X_-\to\Y_N$ is a Moore--Penrose pseudoinverse of the operator $\Q:\Y_N\to\X_-$.
\end{remark}
Let us elaborate equation (\ref{wx}).  The optimal  process $\w x$ can be expressed as  \baaa \w x(t)=\w x
(t,q,s)=\frac{\O}{\pi}\sum_{k\in \ZZ_N }\w y_k \sinc(k\pi+\O t). \label{wxx}\eaaa
Here $\w y=\{\w
y_k\}_{k=-N}^N$ is defined as    \baa \w y=R^{-1}\Q x.\label{wy}\eaa The space $\Y_N$ is finite dimensional, and the operator $R$ can be represented via a matrix
$R=\{R_{km}\}\in \C^{2N+1,2N+1}$, where $k,m=-N,-N+1,...,N-1,N$.  In this
setting, $(Ry)_k=\sum_{k=-N}^N R_{km}y_m$, and  the components of the matrix $R$  are defined as  \baaa R_{km}=\index{
\frac{1}{(2\pi)^2}\sum_{t=q}^s\overline{\left(\frac{e^{im\pi+i\O t}-
e^{-im\pi-i\O t}}{im\pi/\O+i t}\right)} \left( \frac{e^{ik\pi+i\O
t}- e^{-ik\pi-\O t}}{ik\pi/\O+i t}\right) \\=\frac{\O^2}{(2\pi)^2}
\sum_{t=q}^s\frac{e^{-im\pi-i\O t}- e^{imt\O t}}{-im\pi-i\O t}\cdot
\frac{e^{ik\pi+i\O t}- e^{-ik\pi-\O t}}{ik\pi+i\O t}
\\= \frac{\O^2}{(2\pi)^2}\sum_{t=q}^s\frac{-2i\sin(m\pi+\O t)}{-im\pi-i\O t}\cdot
\frac{2i\sin(k\pi+\O t)}{ik\pi+i\O t} \\=
\frac{\O^2}{\pi^2}\sum_{t=q}^s\frac{\sin(m\pi+\O t)}{m\pi+\O t}\cdot
\frac{\sin(k\pi+\O t)}{k\pi+\O t} \\=}
\frac{\O^2}{\pi^2}\sum_{t=q}^s\sinc(m\pi+\O t)\,\sinc(k\pi+\O t)
.\label{R}\eaaa
Respectively, the components of the vector  $\Q^*x=\{(\Q^*x)_k\}_{k=-N}^N$   are defined as
\baa (\Q^*x)_{k}= \index{\sum_{t=q}^s\overline{\left(\frac{e^{ik\pi+i\O
t}- e^{-ik\pi-i\O t}}{ik\pi/\O+i t}\right)} x(t)   =\frac{\O}{\pi}
\sum_{t=q}^s\frac{\sin(k\pi+\O t)}{k\pi+\O t}x(t)
 \\=}\frac{\O}{\pi} \sum_{t=q}^s\sinc(k\pi+\O t)x(t) .
\label{r}\eaa

\begin{remark}
 The process $\w x(t)$  represents the output of a linear causal smoothing filter.  \index{ for the setting  without memory loss and with dynamically expanding  window of observations, i.e., where $q$ is fixed and  $s$ is running.}   It can be noted that the operators
$R$ and $\Q$ have to be recalculated for each $s$, and the values $\w x(t)=\w x_{q,s}(t)$ calculated  for observations $\{x(t),\ s\le t\le q\}$,
 can  be different from the values $\w x_{q,s+\tau}(t)$ calculated for the same $t$ using
 the  observations  $\{x(t),\ s\le t\le q+\tau\}$, where  $\tau>0$. Therefore, this filter
 is not time
invariant.
  \index{ However, if we consider a model  with running  $s$ and with fixed memory window,  i.e., where the size $s-q$ is fixed, this filter will be time invariant. To see this, it suffices to assume that  $s$ is fixed and represent the
running  process
$x(t)$ by the process from $\ell_2(q,s)$ via a time shift.}
\end{remark}
\begin{remark} We have  excluded the case where $\O=\pi$; this case leads to the trivial solution with $x(-t)=y_t$ for $t\in \ZZ_N $.
\end{remark}
\subsection{Detecting left bandlimitness}
Theorem \ref{ThP} allows  to verify if $x\in \ell_2(-\infty,s)$ is left band-limited, i.e.
if  the conditions of Definition \ref{defLBL}  hold. This can be formulated as the following..
\begin{theorem}\label{corr0}
If $x\in\XN$ then $x\in\LBLN$
if and only if
\baaa
(x_L,x_L)_{\ell_2(-\infty,s)}=(\Q x_L,R^{-1}\Q x_L)_{\ell_2}.
\eaaa
\end{theorem}
An alternative approach to detection of left bandlimitness was suggested in \citet{D17}.
 \subsection{Tikhonov regularization}
 This let us consider a modification  of the original optimization problem (\ref{min})
  with penalty on the norm of the solution that  restrains the norm of the solution.  More precisely, let us
 consider the following problem;
 \baa
&&\hbox {Minimize}\quad
\|\w x-x\|^2_{\X_-}+\e\|\w x\|_{\ell_2}^2\quad\breakk\hbox{over}\quad \w x=\XN,\label{minx} \eaa
where $\e>0$ is a parameter.
\begin{theorem}
\label{ThT}
 Problem (\ref{minYe}) has a unique solution \baa
\w x=\Q R_\e^{-1} \Q^*x,\label{wxe}
\eaa
where
\baaa
R_\e=R+\e I,
\eaaa
where $I$ is the unit matrix in $\R^{N\times N}$.
\end{theorem}

 Problem  (\ref{minx}) can be considered as a  regularization of the original optimization problem (\ref{min}) similarly to the setting from \citet{TH}.  We found in numerical experiments that, for large $N$, numerical calculation of inverse matrix $R^{-1}$ is not exact, and
the error $E=\|\Q^*x-R\w y\|_{\ell_2(q,s)}$ for
a numerical solution $\w y$ of equation  (\ref{wy}) does  not vanish.
It appears that the numerical stability can be improved via this regularization.

In fact, the replacement of $R$ by $R_\e$ may help to  decrease the error $E$ even in the original setting.
In particular, we observed that,  for large $N$ and small $\e>0$,
\baaa
\|\Q^*x-R\w y_\e\|_{\X_-}< \|\Q^*x-R\w
y\|_{\X_-}
\eaaa
 for $y_\e=R_\e^{-1}\Q^*x$ and $y=R^{-1}\Q^*x$  such that the optimal solution for problem
(\ref{ThT}) is $\w x=\Q y_\e$ and the optimal solution for problem
(\ref{ThP}) is $\w x=\Q y$. We observed this, for example,  for $N=200$ and $\e=0.05$.

\section{Proofs}
{\em Proof of Lemma \ref{propU}}. Let us prove statement (i).
Let $D\defi\{z\in\C: |z|< 1\}$. Let  $H^2(D)$ be the Hardy space of functions that are holomorphic on
$D$ with finite norm
$\|h\|_{\H^2(D)}=\sup_{\rho<1}\|h(\rho e^{i\o})\|_{L_2(-\pi,\pi)}$.
  Without a loss of generality, we assume that
 $s=0$. In this case,
$\TT=\{t:\ t\le 0\}$.
 It suffices to
prove that if $x\in\XN $ is such that $x(t)=0$ for
$t\le 0$,  then  $x(t)=0$ for $t>0$. Let $X=\Z x$. Since $x\in\XN$, it follows that  $X\in\BNN$.
We have that $X|_D=(\Z x)|_D\in H^2(D)$. Hence, by the property of the Hardy space,
  $X\equiv 0$; see, e.g., Theorem 17.18 from \cite{Rudin}.  This
completes the proof of Lemma \ref{propU} for $N=+\infty$ and
$q=-\infty$. It can be noted that the proof for this case  follows also from the predictability of the two-sided band limited processes
established in \cite{D12a,D12b}.

Let us prove statement (ii). We use an approach
based on non-singularity of special sink matrices established in \citet{LeeF}.  Without a loss of generality, we
assume that $s=N$. \par
Let us consider
first the case when $s-q=2N+1$. It suffices to consider $q=-N$ only; in this
case, the set $\TT=\{t:\ q\le t\le s\}=\{t:\ -N\le t\le N\}$, i.e.,
$\TT=\ZZ_N $ and it has $2N-1$ elements.
 It suffices to
prove that if $x(\cdot)\in\XN $ is such that $x(t)=0$ for
$t\in\TT$, then $x(t)=0$ for $t>0$. By the supposition,  we have that
\baa x(t)=\sum_{k\in \ZZ_N }a_{t,k}y_k =0,\quad  -N\le t\le N,\label{ay0}\eaa
for some set $\{y_k\}$, where $a_{t,k}=\sinc(k\pi+\O
m)$; see, e.g., (\ref{sinc}). By Theorem 1(a) from \citet{LeeF},  the matrix $\{a_{t,k}\}_{k,m=-N}^N\in\R^{2N+1,2N+1}$ is non-singular. Therefore,  linear system
(\ref{ay0}) is a system with a non-singular matrix. Hence $y_k=0$
for all $k$. This proves  Lemma \ref{propU}
for the case where $s-q=2N+1$.

Let us consider the case where $s-q>2N+1$.
In this case, the linear system (\ref{ay0}) has to be considered jointly
with the system \baa\sum_{k\in \ZZ_N }a_{t,k}y_k =0,\quad  -q\le t<-
N.\label{ay00}\eaa Clearly, system (\ref{ay0})-(\ref{ay00}) admits
only a zero solution again. This  completes the proof of Lemma
\ref{propU}. $\Box$

{\em Proof of Corollary \ref{rem1}}.
  Assume first that $s-q=2N+1$. Again, we assume  that $s=N$.
 Since homogeneous
linear system (\ref{ay0}) allows only zero solution, it follows that the
non-homogeneous system \baa\sum_{k\in \ZZ_N }a_{t,k}y_k =x(t_k),\quad
-N\le t\le N\label{ay1}\eaa admits a unique solution $\{y_k\}$ for
any set $\{x(t_k)\}$. Therefore, we proved that
$\{x(t)\}_{t\in\TT}\in \XN $ for any $x\in\ell_2$.  Further, if
$s-q<2N+1$, then there are many solutions of (\ref{ay1}), and these
solutions form a linear manifold.
  This completes the proof of Corollary \ref{rem1}. $\Box$
\par
Consider the mapping $\zeta:\BN \to \XN $ such that
$x(t)=(\zeta (X))(t)=(\Z^{-1} X)(t)$ for $t\in\TT$. It is a linear
continuous operator. By Lemma\ref{propU}, it is a bijection.
\par
 {\em
Proof  of Theorem \ref{Th1}.} The quadratic form here is defined on a finite dimensional
linear subspace of
$\ell_2(q,s)$.  Hence  there exists a unique projection $\w
x$ of $\{x(t)\}_{t\in\TT}$ on $\XN $, and statement (i)  is
proven. Statement
of Theorem \ref{Th1} (ii) follows  from  Lemma \ref{propU}.  $\Box$
\par
{\em Proof of Corollary \ref{rem2}} follows immediately from   Corollary \ref{rem1}
and Lemma \ref{propU}.  $\Box$
\par {\em Proof  of Theorem \ref{ThP}.} Statement (i) follows from
Theorem 1(a) from \citet{LeeF} applied to the matrix $\{\sinc(k\pi+\O
m)\}_{k,m=-N}^N$.  It can be seen from the following. Let $s=N$, $q\le -N$,  $V_0=\{\sinc(k\pi+\O
m)\}_{k=-N,m=q}^{N,s}\in\R^{2N+1,q-s}$, then $R-V_0V_0^\top$ is non-negative definite. It follows that $R$ is positively defined. This proves statement (i).
\par
Let us prove statement (ii).
Let the Hermitian form $F:\XN \times
\X_-\to\R$ be defined as\baaa F(\w x,x)=\|\w x-x\|_{|X_-}^2=\sum_{t=q}^s|\w
x(t)-x(t)|^2 . \eaaa
Further, let the Hermitian form $G:\Y_{N}\times \X_-\to\R$ be defined as \baaa
G(y,x)=F(\Q y,x).
\label{formG}\eaaa It follows that \baaa
G(y,x)=\|\w x-x\|_{\X_-}^2 , \quad \w x=\Q y.
\label{formGG}\eaaa
Clearly, problem (\ref{min}) can be replaced by the minimization problem \baaa
\hbox{Minimize}\quad && G(y,x)\quad\hbox{over}\quad y\in
\Y_{N}.\label{minY} \eaaa
By the definition, it follows that  \baaa &&G(y,x)=(\Q y-x,\Q y-x)_{\X_-}\breakk= (\Q y,\Q y)_{\X_-}-2\Re(\Q y,x)_{\X_-}+(
x,x)_{\X_-}\nonumber\\&&= (\Q y,\Q y)_{\X_-}-2\Re(\Q y,x)_{\X_-}+(
x,x)_{\X_-}\hphantom{xxx}.
\eaaa
As was mentioned above, it follows from the properties of the sinc function that the mapping $\Q^*:\ell_2\to\ell_2$ defined by (\ref{Q*}) is  continuous, and, therefore,
the operator $\Q^*$ can be extended as a continuous linear operator $\Q^*:\X_-\to \Y_N$.
 It follows that  \baaa  G(y,x) =(y,Ry)_{\Y_N}-2\Re(y,\Q^*x)_{\Y_N}+(
x,x)_{\X_-},\label{G0}\hphantom{xxx}\eaaa
i.e., this is a quadratic form defined on $\Y_N\times\Y_N$. By the definitions, the operator $R$ is non-negative definite, and, by Lemma \ref{propU},
\baa
(y,Ry)_{\Y_N}=(\Q y,\Q y)_{\X_-}>0\quad \forall y\neq 0_{\Y_N}.\label{G}
\eaa Finally, statement (ii) follows from the invertibility of $R$ and the standard properties of the quadratic forms.  This completes the proof of Theorem
\ref{ThP}. $\Box$
\par
{\em Proof of Theorem \ref{corr0}.} By Theorem \ref{ThP} applied for $q=-\infty$, $s=0$, 
$x_L$ is an unique solution of problem (\ref{min}). In this case, $\X_-=\ell_2(-\infty,s)$. By (\ref{G}), $x_L$ is an unique solution of problem
\baaa
 \hbox{Minimize}\quad
  \|x-\Q y\|^2_{\ell_2(-\infty,s)}\quad\hbox{over}\quad y\in\ell_2.
  \label{optL}\eaaa
By the definitions,
 \baaa
 &&\|x-\Q y\|^2_{\ell_2(-\infty,s)}\breakk=(y,Ry)_{\ell_2} -2(y,\Q x_L)_{\ell_2} +(x_L,x_L)_{\ell_2(-\infty,s)}.
  \eaaa The value \baaa
 \rho=(x_L,x_L)_{\ell_2(-\infty,s)}-(\Q x_L,R^{-1}\Q x_L)_{\ell_2}\eaaa
 represents the optimal value of problem (\ref{optL}). Clearly, $x_L$ is left band-limited if and only if
 $\rho=0$. This completes the proof. $\Box$
 \par
{\em Proof of Theorem \ref{ThT}}. By the definition of $\Q$, for $\w X=\Q y$ and $\w x=\Z^{-1}\w X$, we have that
\baaa
\|\w x\|_{\ell_2}^2=\frac{1}{2\pi}\|\w X\|_{L_2(-\pi,\pi)}^2=\|y\|_{\ell_2}^2.
\eaaa
Hence problem (\ref{minx}) is equivalent to the problem
\baa
&&\hbox {Minimize}\quad
\|\w x-x\|^2_{\X_-}+\e\|y\|_{\ell_2}^2\quad\breakk\hbox{over}\quad y\in
\Y_{N}, \quad \w x=\Q y.\label{minYe} \eaa
The remaining part of the proof is similar to the proof of Theorem \ref{ThP} with $R$ and $G(y,x)$
replaced by $R_\e$ and $G(y,x)+\e\|y\|_{\ell_2}^2$ respectively. $\Box$
\section{Some numerical experiments}
We did some  numerical  experiments to compare statistically the performance of our band-limited  extrapolations with
extrapolations based on splines applied to causally smoothed processes.
\subsection{Simulation of the input processes}
       The setting of Theorems \ref{ThP} does not involve stochastic processes and probability measure; it is oriented on extrapolation of real sequences.
      However, to  provide  sufficiently large  sets of input sequences for statistical estimation,  we used processes $x$ generated via Monte-Carlo simulation
      as a stochastic  process evolving as
       \def\b{\beta}
      \baaa
       &&z(t)={\rm A}(t) z(t-1)+\s\eta(t), \quad  t\in\ZZ,\qquad \breakk  x(t)=c^\top z(t).
       \label{AR2}\eaaa
       Here $z(t)$ is a process with the values in $\R^\nu$, where $\nu\ge 1$ is an integer, $c\in \R^\nu$.
       The process  $\eta$ represents a noise  with values in $\R^\nu$,
          ${\rm A}(t)$ is a matrix with the values in
        $\R^{\nu\times \nu}$ with the spectrum inside $\T$, $\s>0$.
        The matrices ${\rm A}(t)$ are switching  values randomly at random times; this replicates a  situation where the parameters of a system
        cannot be recovered from the observations such as described in the review \cite{SFS}.

In each simulation, we selected random and mutually independent $(\nu,\s, z(-N),A(\cdot),\eta)$, with vectors and
 matrices having  mutually independent components. We selected $\nu\in\{1,...,10\}$ randomly with equal probability, and we selected $\s$ from the uniform distribution
 on the interval $(0,2)$.
   The process $\eta$ was selected as a stochastic discrete time Gaussian white noise with the values in $\R^\nu$
such that $\E \eta(t)=0$ and $\E |\eta(t)|^2=1$. The initial vector  $z(-N)$ was selected randomly  with the components
     from the uniform distribution on $(0,1)$.
   The components of the matrix ${\rm A}(-N)$ was selected from the uniform distribution on $(0,1/\nu)$. Further, to simulate randomly changing ${\rm A}(t)$,
     a random variable $\xi$ distributed
     uniformly on $(0,1)$ and  independent on $({\rm A}(s)|_{s<t},\eta,z(-N))$ was simulated for each time $t>-N$.
         In the case where $\xi<0.5$,  we selected
     ${\rm A}(t)={\rm A}(t-1)$. In the case where $\xi\ge 0.5$,
     ${\rm A}(t)$ was simulated randomly  from the same distribution as ${\rm A}(-N)$, independently on $({\rm A}(s)|_{s<t},\eta,z(-N))$.
This setting  with randomly changing  ${\rm A}(t)$ makes impossible to identify the parameters of equation (\ref{AR2})  from the current observations.

\subsection{Comparison of band-limited extrapolation with spline extrapolations}
We compared root-mean-square errors (RMSEs)  for the forecasting  via extrapolations of the band-limited approximation
obtained in Theorem \ref{ThP} with the RMSEs of standard some spline  extrapolations.

\def\EE{\mathbb{E}}
We denote below by $\EE$ the sample mean across the Monte Carlo trials.

We estimate the root-mean-square error (RMSE) for the forecasting
    \baa
      e_\BL=\EE\left[\left(\sum_{t=q+1}^{q+L}|x(t)-\w x_\BL(t)|^2\right)^{1/2}\right],\label{err}\eaa
given that the extrapolation of the band-limited approximation is accepted as
the forecast. Here  $\w x_\BL$ is an extrapolation on future times $t=q+1,q+2,...,q+L$ of the band-limited approximation described in Theorem
\ref{ThP} for the underlying process $x|_{t=q,...,s}$  with $R$ replaced by $R_\e=R+\e I$, in the terms of this theorem.
 The choice of integers $L>0$ defines the forecasting horizon.

 We have compared these values with similar values obtained for some standard spline extrapolations
 of the causal $h$-step moving average process for $x$. More precisely, to
 take into the account truncation, we used a modification of the causal moving average \baaa
  x_\MA(t)=\frac{1}{\min(h,t+N+1)}\sum_{k=\max(t-h,-N)}^t x(k),\quad t\ge -N.
 \eaaa
 \par
 For three different types of standard spline extrapolations, we calculated the root-mean-square error (RMSE)
 \baa
      e_\spl=\EE\left[\left(\sum_{t=1}^L|x(t)-x_\spl(t)|^2\right)^{1/2}\right], \label{err2}\eaa
 given that the spline is accepted as
the forecast. Here $x_\spl$ is a spline extrapolation  of the observed moving average $x_\MA|_{t=q,...,s}$.
 We considered the
 the  piecewise cubic extrapolation,
 the shape-preserving piecewise cubic extrapolation, and
 the linear  extrapolation.   We used built in MATLAB code {\em interp1} for calculation of these extrapolations.

 We used  smoothed
  moving average process $x_\MA$ as inputs because we found that applications directly to the original
 "noisy" process $x$ produces  quite unsustainable extrapolation with  large errors $e_\spl$.
\par
We calculated and compared $e_\BS$ and $ e_\spl$.
Table \ref{simulation} shows the ratios  $e_{\BL}/e_\spl$ for some combinations of parameters. For these calculations, we used  $c=(1/\nu,1/\nu,...,1/\nu)^\top$, $h=10$, and $\e=0.1$.

\subsection{Impact of preliminary smoothing}
It is common to apply a forecasting method to  processes that are preliminary smoothed by a causal filter.
For many methods, it helps to improve performance.
We did  some experiments to investigate the impact of this smoothing on relative performance of the band-limited extrapolation and
spline extrapolation. We repeated  experiments  described above with the following modification: we calculated band-limited projections
and their extrapolations for the  causally smoothed process $x_\MA(t)$, and compared the performance of the corresponding predictor
  with the performance of defined by the same spline extrapolations and applied to $x_\MA$ as described above.
  Again, we calculated and compared corresponding RMSEs
  $e_\BL$ and $ e_\spl$, with  $c=(1/\nu,1/\nu,...,1/\nu)^\top$, $h=10$, $s=0$, $q=-150,-70$, and $\e=0.1$.
Table \ref{simulation2} shows the ratios  $e_{\BL}/e_\spl$ for some combinations of parameters with $s=0$ and $q=-150, -70$.

\par
Comparing Tables \ref{simulation} and \ref{simulation2}, we observe that
using smoothed moving average process $x_\MA$ instead of the original process $x$ as the input for the band-limited
approximation and extrapolation leads to slightly increased ratios $e_\BL/e_\spl$.  Since the errors $ e_\spl$ used for
Table \ref{simulation2} are the same as the ones used for Table \ref{simulation}, it follows  that the errors $e_\BL$ are slightly larger with the
smoothed moving average process $x_\MA$ is used instead of the original process $x$ as the input for the left band-limited
approximation and extrapolation.  The fact that a preliminary smoothing does not improve performance of the suggested left band-limited
 smoothing  speaks  in favor of our method.  
This  is expected for filters targeting approximation of an ideal low-pass filter, because
moving average damps higher frequency but distorts significantly  a signal characteristics on a wider spectrum.

\begin{table}[h]
\caption{The ratios of root-mean-square errors $e_{\BL}/ e_\spl$ defined by (\ref{err})-(\ref{err2}) for band-limited extrapolation
of  $x|_{q\le t\le 0}$  and spline extrapolations of  $x_\MA|_{q\le\t\le 0}$.}
  \centering
\begin{tabular}{c|cccc}
  \hline &\multicolumn{4}{c}{$e_\BL/e_\spl$} \\
  \hline
   &\multicolumn{4}{c}{Extrapolation horizons }\\

$\begin{array}{c}
   \hbox{Type of spline}\vspace{-1mm} \\
   \hbox{ extrapolation}
 \end{array}$   &$L=1$&$L=4$&$L=8$&$L=12$\\
  \hline \\&\multicolumn{4}{c}{$\O=\pi/4$, $N=50$, $q=-150$}
  \\
Linear   & 0.8900 & 0.8825 & 0.7781  & 0.6742\\
Piecewise cubic  & 0.8429 & 0.2359 & 0.052   & 0.0192\\
$\begin{array}{c}
   \hbox{Shape-preserving}\vspace{-1mm} \\
   \hbox{piecewise cubic}
 \end{array}$ & 0.9001 & 0.6303 & 0.1414   & 0.0448 \\
\\
&\multicolumn{4}{c}{  $\O=\pi/2$, $N=30$, $q=-70$ } \\
Linear   & 0.9010 & 09005. & 0.7840  & 0.6837\\
Piecewise cubic  & 0.8675 & 0.2408 & 0.0530   &0.0193 \\
$\begin{array}{c}
   \hbox{Shape-preserving}\vspace{-1mm} \\
   \hbox{piecewise cubic}
 \end{array}$ & 0.9082 & 0.6475 & 0.1436   &0.0464
\end{tabular}
\label{simulation}
\end{table}

\comm{LINEAR [vec]=FstatisSS(0.25*pi,50,50,12,0.1,1,30000,12)vec =0.8900    0.8825    0.7781    0.6742

[vec]=FstatisSS(0.5*pi,30,10,12,0.1,1,30000,12)vec = 0.9010    0.9005    0.7840    0.6837

SPLINE

[vec]=FstatisSS(0.25*pi,50,50,12,0.1,1,30000,12)vec =0.8429    0.2359    0.0525    0.0192

[vec]=FstatisSS(0.5*pi,30,10,12,0.1,1,30000,12)vec =0.8675    0.2408    0.0530    0.0193

PCHIP

[vec]=FstatisSS(0.25*pi,50,50,12,0.1,1,30000,12)vec =0.9001    0.6303    0.1414    0.0458

[vec]=FstatisSS(0.5*pi,30,10,12,0.1,1,30000,12)vec =0.9082    0.6475    0.1436    0.0464
}
\index{\vspace{2mm}
{\sm $\EE e_d$ represent the average distances  (\ref{err2}) from the future process:  $\EE e_1$ for
 the  piecewise cubic extrapolation $\ww x_2$ of the moving average, $\EE e_1$ for the shape-preserving piecewise cubic extrapolation,
 $\EE e_3$ for the linear extrapolation; $L$ is the extrapolation horizon.}}

\begin{table}[h]
\caption{The ratios of root-mean-square errors $e_{\BL}/ e_\spl$ defined by (\ref{err})- (\ref{err2}) for band-limited extrapolation
of the moving average $x_\MA|_{q\le t\le 0}$  and spline extrapolations of  $x_\MA|_{q\le\t\le 0}$.}
  \centering
\begin{tabular}{c|cccc}
  \hline &\multicolumn{4}{c}{$e_\BL/e_\spl$} \\
  \hline
   &\multicolumn{4}{c}{Extrapolation horizons }\\
$\begin{array}{c}
   \hbox{Type of spline}\vspace{-1mm} \\
   \hbox{ extrapolation}
 \end{array}$   &$L=1$&$L=4$&$L=8$&$L=12$\\
  \hline \\&\multicolumn{4}{c}{$\O=\pi/4$, $N=50$, $q=-150$}
   \\
Linear     & 0.9891 & 0.8915 & 0.7719  & 0.6677\\
Piecewise cubic  & 0.9464 & 0.2425 & 0.0536   &0.0196 \\
$\begin{array}{c}
   \hbox{Shape-preserving}\vspace{-1mm} \\
   \hbox{piecewise cubic}
 \end{array}$ & 0.9882 & 0.6490 & 0.1446   & 0.0468 \\
&\multicolumn{4}{c}{ $\O=\pi/2$, $N=30$, $q=-70$} \\
Linear   & 0.9736 &  0.8934 & 0.7785 & 0.6777\\
Piecewise cubic  & 0.9362 & 0.2416 & 0.0530   & 0.0193\\
$\begin{array}{c}
   \hbox{Shape-preserving}\vspace{-1mm} \\
   \hbox{piecewise cubic}
 \end{array}$ & 0.9783 & 0.6531 & 0.1459   &0.0472 \\
\end{tabular}\label{simulation2}
\end{table}

\comm{LINEAR[vec]=FstatisSS(0.25*pi,50,50,12,0.1,1,30000,12)vec =0.9891 + 2.6530i   0.8915 + 1.2824i   0.7719 + 0.7295i   0.6677 + 0.4839i

[vec]=FstatisSS(0.5*pi,30,10,12,0.1,1,30000,12)vec =0.9736 + 2.6398i   0.8934 + 1.3931i   0.7785 + 0.7471i   0.6777 + 0.4958i

SPLINE:[vec]=FstatisSS(0.25*pi,50,50,12,0.1,1,30000,12) vec =0.9464 +i...   0.2425 + i...   0.0536 + i...   0.0196 + i...

[vec]=FstatisSS(0.5*pi,30,10,12,0.1,1,30000,12)vec =0.9362 + i...    0.2416 + i...    0.0530 + i...    0.0193 + i...

PCHIP;

[vec]=FstatisSS(0.25*pi,50,50,12,0.1,1,30000,12)vec =0.9882    0.6490    0.1446    0.0468

[vec]=FstatisSS(0.5*pi,30,10,12,0.1,1,30000,12)vec =0.9783    0.6531    0.1459    0.0472
}

\subsection{Discussion of the results of the experiments}
 The experiments demonstrated a good numerical stability of the  method;
 the results were quite robust with respect to truncation of the input processes and deviations of parameters.
 For each  entry in Tables \ref{simulation}-\ref{simulation2}, we used  30,000
     Monte-Carlo trials; increasing or decreasing  the number of Monte-Carlo trials also gives very close results.
 For instance,
an experiment with 60,000 Monte-Carlo trials produced the set of results  $(0.9834,    0.6436,    0.1426,    0.0459)$ for the last row of
Table \ref{simulation2}. An experiment with 15,000 Monte-Carlo trials produced the set of results (0.8508,0.2378, 0.0530,0.0194) for the 3rd row of
Table \ref{simulation}.

The ratios   $ e_{\BL}/ e_\spl$ are decreasing further as the horizon $L$ is increasing, hence we omitted the results for $L>12$.
By the same reasons, we omitted results with classical extrapolations applied directly to $x(t)$ instead of the moving average $x_\MA(t)$, since
the distance (\ref{err2}) is quite large in this case due the presence of the noise.

\begin{figure}[ht]
\centerline{\psfig{figure=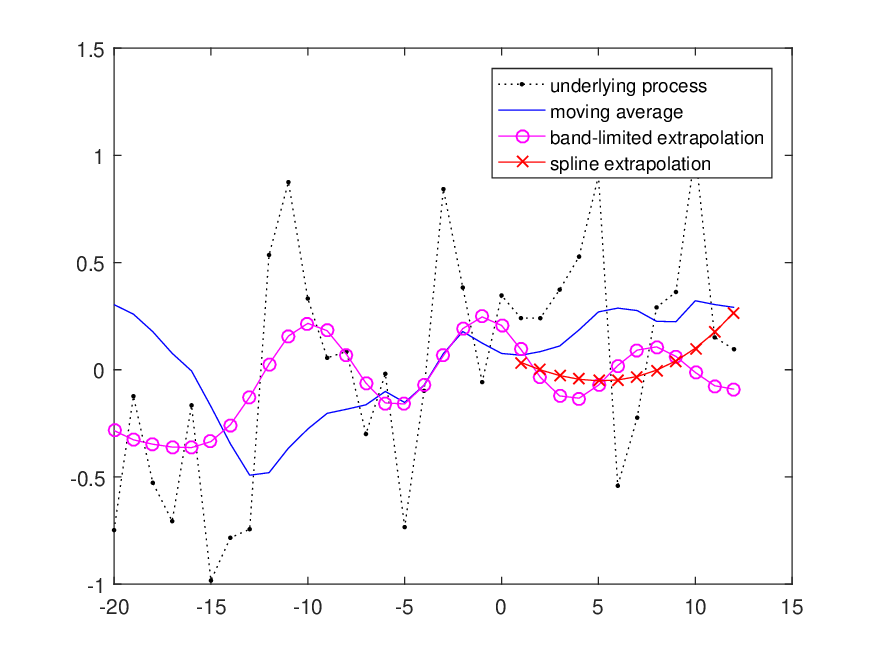,width=7.5cm,height=6.0cm}}
\centerline{\psfig{figure=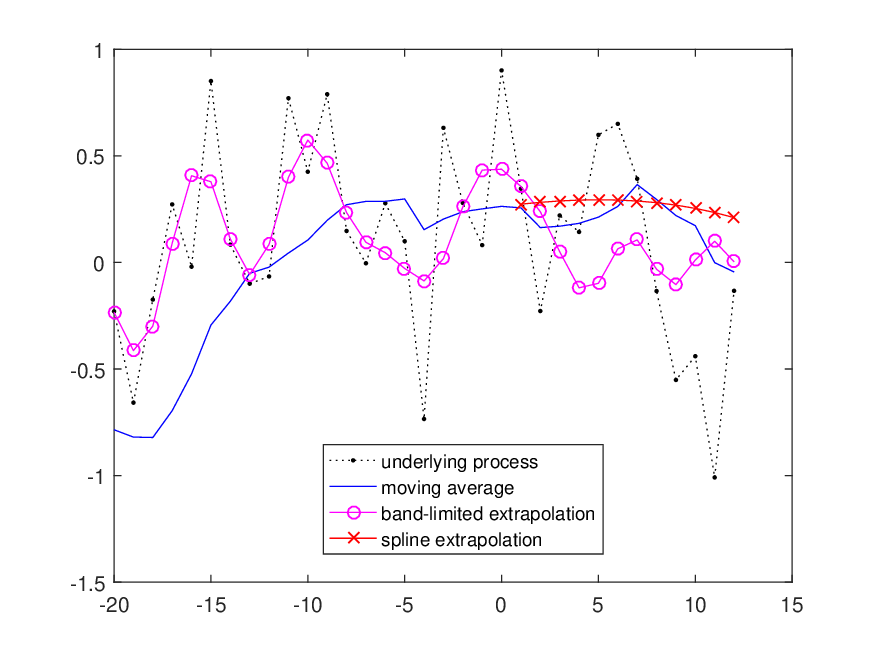,width=7.5cm,height=6.0cm}}
\caption[]{\sm Example of a path $x(t)$, its moving  average $x_\MA$,
the band-limited causal approximation and extrapolation $\w x_\BL$
of  $x|_{t\le 0}$,
and shape-preserving piecewise cubic extrapolation $x_\spl$ of $x_\MA|_{t\le 0}$,
 with
$\O=\pi/4$, $N=50$,  $q=-150$, $s=0$, $h=10$, $\e=0.1$ (top) and
$\O=\pi/2$, $N=30$, $q=-70$, $s=0$, $h=10$, and $\e=0.1$ (bottom).  } \vspace{0cm}\label{fig-1}
\end{figure}

Figure  \ref{fig-1} shows examples of paths $x(t)$, their band-limited causal approximation and extrapolation $\w x_\BL$,
with the same parameters as were used for Tables \ref{simulation}.
These figures also show
the moving  averages $x_\MA(t)$, and their spline extrapolations  $x_\spl(t)$.

\begin{figure}[ht]
\centerline{\psfig{figure=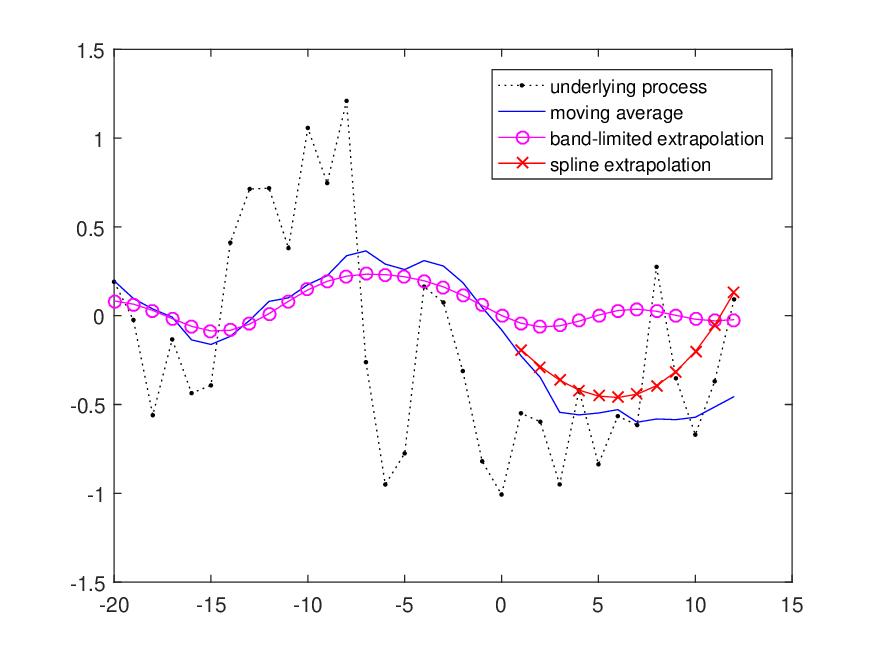,width=7.5cm,height=6.0cm}}
\centerline{\psfig{figure=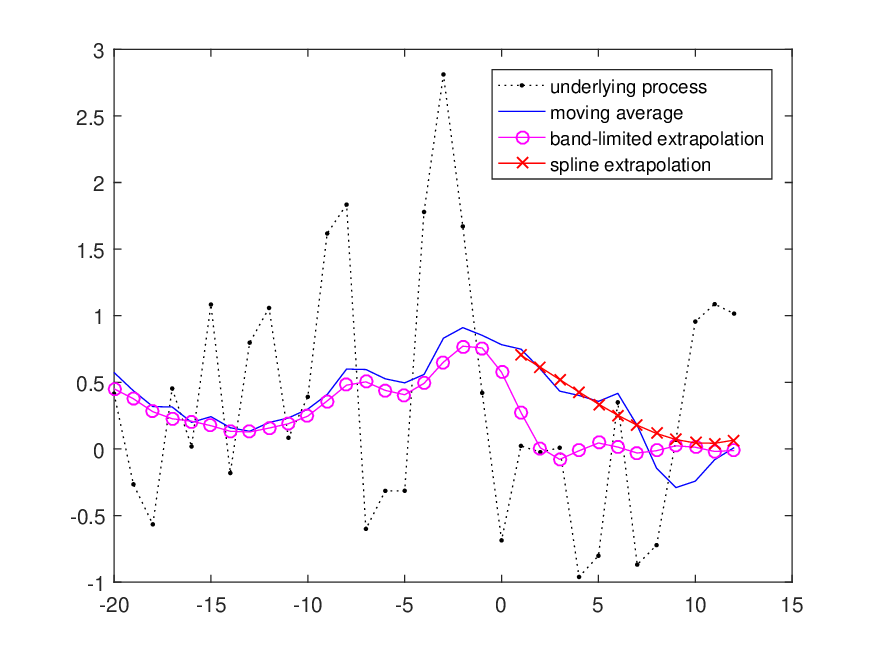,width=7.5cm,height=6.0cm}}
\caption[]{\sm Example of a path $x(t)$,  its moving  averages $x_\MA(t)$,
the  band-limited causal approximation and extrapolation $\w x_\BL$ of $x_\MA|_{t\le 0}$,
 and shape-preserving piecewise cubic extrapolation $x_\spl$ of $x_\MA|_{t\le 0}$ with
$\O=\pi/4$, $N=50$,  $q=-150$, $s=0$, $h=10$, $\e=0.1$ (top) and
$\O=\pi/2$, $N=30$, $q=-70$, $s=0$, $h=10$, and $\e=0.1$ (bottom). } \vspace{0cm}\label{fig-2}
\end{figure}

Figure  \ref{fig-2} shows examples of paths $x(t)$, the moving  averages $x_\MA(t)$,
the  band-limited causal approximation and extrapolation $\w x_\BL$ obtained for $x_\MA(t)$,
with the same parameters as were used for Tables \ref{simulation2}.
The figure also shows spline extrapolations  $x_\spl(t)$ of $x_\MA$.

\section{Possible applications and future development}
The  approach suggested in this paper  allows many modifications. We
outline below some possible straightforward modifications as well as
  more challenging problems and possible applications  that we leave for the future research.
\begin{enumerate}
\item
The mean-square optimal causal band-limited approximations of one-sided
sequences suggested above  can be interpreted as outputs causal smoothing filter.
To accommodate the current
flow of observations, the coefficients of the sinc series have to
be changed dynamically; therefore,  the corresponding filter is not time
invariant.   It can be noted that, for some problems, time
invariance for a filter is not crucial. For example, a typical
approach to forecasting in finance is to approximate the known path
of the stock price process by a process that has a unique
extrapolation that can be used as a forecast. This
procedure can be done at current time; it is not required that
the same forecasting rule will be applied at future times.

\item
Tables \ref{simulation}-\ref{simulation2}  show that the band-limited extrapolation performs better than the spline extrapolations;
some additional experiments with other choices of parameters demonstrated the same trend.
However, experiments did not involve more advanced methods beyond the listed above spline methods.
Nevertheless, regardless of the results of these experiments,
 potential importance of band-limited extrapolation is self-evident  because  its
physical meaning:  a band-limited part can be considered
 as a regular part of a process purified from a noise represented by high-frequency component. This is controlled
 by the choice of the band. On the other hand, the choice of particular splines does not have a physical
 interpretation.

\item
The set $\{e^{i\o}, \ \o\in [-\O,\O]\}\subset\T$ can be replaced by another set, for example, by a set that is not
necessarily connected, in a setting that is close to one from \citet{F94}. This would require a minor modification of the algorithm.
\item It is possible to consider a setting where some observations of the past values  $x(t)$  are missing.
\item
Instead of Fourier series, expansion by another
basis in $L_2(-\O,\O)$  can be used, for instance, such as suggested in \citet{TH}.
 The space
$L_2(-\O,\O)$ can be replaced by a weighted
$L_2$-spaces, for a  weight  representing a relative importance of the
approximation on different frequencies.

\end{enumerate}
\subsection*{Acknowledgements}
This work  was supported by ARC grant of Australia DP120100928 to the author.


\end{document}